\documentclass[12pt]{article}
\usepackage{amsfonts}
\usepackage{amssymb}
\usepackage{amsmath}
\usepackage{url}

\begin{document}

{\footnotesize
\noindent
This is author's translation  of his paper
originally published in Russian in

\emph{Diskretnyi Analiz}, issue 26,  65-71 (1974),

Novosibirsk, Institute of Mathematics

of Siberian Branch of
Academy of Sciences of the USSR;

\url{http://www.zentralblatt-math.org/zmath/en/advanced/?q=an:0298.05114}

\noindent
Russian original can be downloaded from

\url{http://logic.pdmi.ras.ru/~yumat/papers/22_paper/22_home.html}

\noindent
Part of the results was rediscovered  by
N.\,Alon and M.\,Tarsi in  \emph{Combinatorica},  12:2, pp.\,125--134   (1992);
\url{http://www.ams.org/mathscinet-getitem?mr=1179249}

\noindent
The author is very grateful to Martin Davis
for some help with the English.
}

\

\hfill UDK 519.1

\

\begin{center}
A CRITERION FOR VERTEX COLORABILITY OF A GRAPH\\
STATED IN TERMS OF EDGE ORIENTATIONS\\[5mm]
Yu.V.Matiyasevich
\end{center}

\

1. L.\,M.\,Vitaver [1] and G.\,I.\,Minty [2]
suggested criteria for vertex colorability of a graph in at most a
given number $k$ of colors; these criteria  are stated in terms of the orientation
of the edges.
(Both results are reproduced in the monograph [3] from which we
borrow terminology and notation). One additional criterion of this kind is
given below.

Let us consider all possible directed graphs which can be obtained from a
(non-directed) graph $L=(X,U)$ by orienting each of its edges. For
each such graph $ \overrightarrow{L}$ there is a function
$s^+_{\overrightarrow{L}}(x)$ defined on $X$ whose value is equal to the
number of outgoing edges from the vertex $x$. For two such graphs $ \overrightarrow{L'}$
and $ \overrightarrow{L''}$ we say that they \emph{are congruent modulo $k$}
if for each vertex $x$
\begin{equation*}
s^+_{\overrightarrow{L'}}(x)\equiv s^+_{\overrightarrow{L''}}(x)\pmod{k}.
\end{equation*}
Clearly, this relation is reflexive, symmetric, and transitive; hence all
orientations of the graph $L$ split into \emph{equivalence classes modulo $k$ }.

Let us introduce one additional equivalence relation, namely, let us say
that $ \overrightarrow{L'}$
and $ \overrightarrow{L''}$ \emph{agree} if one of these graphs can be obtained from
the other by changing the orientation of an even number of edges. The relation of
agreement splits every equivalence class modulo $k$ into two subclasses which will be called
\emph{adjacent}.

The proposed criterion is stated below in the form of separate sufficient
and necessary conditions. The necessary condition is formally stronger
than  the sufficient one, so any intermediate condition can serve by itself
as a criterion.

\

For a (non-directed) graph to have a vertex coloring in at most $k$ colors, it is
\begin{itemize}
\item SUFFICIENT that there exist an adjacent pair of subclasses modulo $k$ which disagree with respect to their number of possible orientations modulo $k$;
\item NECESSARY that for every natural number $l$ different from $1$ and
co-prime with $k$ there are two adjacent subclasses modulo $k$ whose
cardinalities are distinct modulo $l$.
\end{itemize}

An interesting intermediate criterion can be obtained for an odd $k$ by
taking $l=2$ (this criterion can be stated without introducing the notion of
adjacent subclass):

\

For a (non-directed) graph to have a vertex coloring in at most $k$ colors,
with an odd $k$, it is necessary and sufficient that at least one of its
modulo $k$ equivalence classes contains an odd number of orientations.

\

The requirement of co-primality of $l$ and $k$ is essential:
a simple circle of an even length has a coloring in 2 colors but each of
its non-empty equivalence classes modulo $2$ consists of 2 orientations.

Let us emphasize the following property of the proposed criteria distinguishing them from
those of Vitaver and Minty. In their criteria, the existence of a coloring
is connected to the existence of another object, namely, an orientation
of a special kind. The coloring and the orientation have a close relationship,
so the graph has few or many colorings corresponding to the existence of
few or many such orientations. In our criteria the existence of a coloring is
also connected with the existence of another object--a pair of adjacent
subclasses with non-equal cardinalities, but there is no close relationship
between such pairs and colorings. The empty graph with $n$ vertices has
$k^n$ colorings, that is, the maximal possible number of colorings in $k$ colors,
but it has only one orientation and hence only one pair of adjacent subclasses
satisfies the criteria, while a graph with $n$ vertices  could have  up to
$k^n$ pairs of adjacent subclasses. On the other hand it can be shown that if
a graph has a unique (up to renaming) coloring, then at least $k^{n-k}$
adjacent subclasses meet the sufficient condition of our criteria. Thus, our
criteria are more
efficient on graphs with few colorings, that is in the cases which usually
are of greatest interest and of greatest difficulty.

Let us mention that the proposed criteria is valid also for graphs with loops
provided that we assume that each loop can be oriented in two ways;
the criteria is also valid for graphs with multiple edges. The proof for the general
case differs only in a slight complication of notation.

\vspace{3mm}
2. We now introduce the notions and notation required for the proof of the criteria.

Let us introduce a one-to-one correspondence between the vertices of graph $L$ and the
formal variables $x_1,\dots,x_n$ (below we just identify the vertices with these variables).
Let us fix an orientation $ \overrightarrow{L}^*=(X^*,\overrightarrow{U}^*)$ of
all of the edges and let $M_L(x_1,\dots,x_n)$ denote the \emph{characteristic polynomial} of
graph $L$ defined as:
\begin{equation}
\prod_{\overrightarrow{x_ix_j}\in \overrightarrow{U}^*}(x_i-x_j)\label{eq1}.
\end{equation}
(Our notation doesn't reflect  the choice of the orientation $ \overrightarrow{L}^*$
but it is easy to see that polynomials corresponding to different choices of
$ \overrightarrow{L}^*$ differ only in sign, and this difference is inessential
in what follows.) If we treat colors as elements of some ring with no divisors of zero
and take the value of a variable $x$ to be equal to the color of vertex $x$, then
the inequality
\begin{equation}
M_L(x_1,\dots,x_n)\neq 0
\label{eq2}
\end{equation}
distinguishes colorings among all the possible ways to assign elements of the ring
to the variables.

Let us temporarily suppose that $q=k+1$ is a prime number. Let us choose as colors the
\emph{non-zero} elements of the finite field $\mathrm{GF}(q)$ with $q$ elements
(that is, the field of residues modulo $q$); however, we will permit the variables
to assume \emph{arbitrary} values from this field. Now the role of the inequality
\eqref{eq2} is played by the inequality
\begin{equation}
x_1\cdots x_n M_L(x_1,\dots,x_n)\neq 0.
\label{eq3}\end{equation}
In other words, the graph $L$ has no coloring in $k$ or fewer colors if and
only if the polynomial
\begin{equation}
x_1\cdots x_n M_L(x_1,\dots,x_n)
\label{eq4}\end{equation}
is \emph{identically equal} to zero.

In a finite field a polynomial can be \emph{identically equal} to zero
without being \emph{formally equal} to the zero polynomial, that is
to the polynomial with all coefficients equal to the zero element of the field.

An example of such a polynomial is given by $x^q-x$ using  Fermat's Little Theorem
according to which
\begin{equation}
x^q\equiv x \pmod{q}.
\label{eq5}\end{equation}
However, a polynomial having degree (in each of the variables)
at most $q-1$ can be identically equal to zero only it is the formally zero polynomial
(for polynomials in one variable this follows from the fact that the number
of roots of a polynomial isn't greater than its degree; this can be easily
generalized to polynomials in many variables by induction on their number).
Below, in order to denote formal (coefficient-wise) equality of two
polynomials we'll use the symbol $\eqcirc $.

Let $i$ and $j$ be two natural numbers such that $i<j$. \emph{Maximal reduction
according to scheme}
$x^j\rightarrow x^i$ of a polynomial
\begin{equation*}
A(x_1,\dots,x_n)\eqcirc
\sum_{j_1,\dots,j_n}a_{j_1,\dots,j_n}x_1^{j_1}\dots x_n^{j_n}
\end{equation*}
is defined as the polynomial
\begin{equation*}
\sum_{j_1,\dots,j_n}a_{j_1,\dots,j_n}x_1^{i_1}\dots x_n^{i_n}
\end{equation*}
where $i_m$ is the least integer which isn't less than $i$ and which is
congruent to $j_m$ modulo $j-i$; the latter polynomial will be denoted
$R_i^j[ A(x_1,\dots,x_n)]$. According to \eqref{eq5}
for every polynomial $A(x_1,\dots,x_n)$   we have the following identity in
the field    $\mathrm{GF}(q)$:
\begin{equation*}
R_q^1[ A(x_1,\dots,x_n)]= A(x_1,\dots,x_n).
\end{equation*}
 On the other hand,  the maximal reduction
according to scheme
$x^q\rightarrow x^1$ has the degree, in each variable, at most $q-1$
and hence the polynomial \eqref{eq4} is identically equal to zero
if and only if the polynomial
\begin{equation*}
R_q^1[x_1\cdots x_n M_L(x_1,\dots,x_n)]
\end{equation*}
is formally equal to the zero polynomial in the  field $\mathrm{GF}(q)$.
It is easy to see that
\begin{equation*}
R_1^q[x_1\cdots x_n M_L(x_1,\dots,x_n)]\eqcirc
x_1\cdots x_n R^{q-1}_0[ M_L(x_1,\dots,x_n)].
\end{equation*}
Multiplication by $x_1\cdots x_n$ transforms a formally zero (non-zero) polynomial
into a formally zero (respectively, non-zero) polynomial. Thus
the polynomial \eqref{eq4} is identically equal to zero
if and only if the polynomial
\begin{equation}
R^{q-1}_0[x_1\cdots x_n M_L(x_1,\dots,x_n)]
\label{eq6}\end{equation}
is formally equal to the zero polynomial in the  field $\mathrm{GF}(q)$.

Let us weaken the assumptions that $q=k+1$ and $q$ is a prime. From now on we assume that
$q$ can be any natural number meeting the following two conditions:
$q\equiv 1 \pmod{k}$; there exists a finite field $\mathrm{GF}(q)$ with
$q$ elements (that is, $q$ should be a power of a prime number). We still allow
variables to take arbitrary values from the field $\mathrm{GF}(q)$
but for colors we now take only those non-zero elements that can be represented in the
form $x^m$ where $m$ is defined from the equality $q=mk+1$. Let us verify that
indeed we have exactly $k$ colors.

Suppose that we had $k'$ colors $a_1,\dots,a_{k'}$. Then each of the $mk$
non-zero elements of the
field $\mathrm{GF}(q)$ is a root of one of the $k'$ equations
\begin{equation*}
x^m=a_1,\qquad\dots,\qquad x^m=a_{k'}.
\end{equation*}
Each such equation has at most $m$ roots, hence $k'\ge k$. On the other hand,
in the  field $\mathrm{GF}(q)$ we have the identity
$$ x^q=x$$ which is a counterpart of \eqref{eq5}. This implies that $a_1,\dots,a_{k'}$ are roots of the
equation $$ x^k=1,$$ hence $k'\le k$. We have established that $k'= k$.

Let us assume that a vertex $x$ is colored by color $x^m$. Now the role of the inequality
 \eqref{eq3} will be played by the inequality
\begin{equation*}
x_1\cdots x_n M_L(x_1^m,\dots,x_n^m)\neq 0,
\end{equation*}
and the role of polynomial \eqref{eq6} will be played by the polynomial
\begin{equation*}
R^{q-1}_0[ M_L(x_1^m,\dots,x_n^m)].
\end{equation*}
It is easy to verify that
\begin{eqnarray*}
R^{q-1}_0[ M_L(x_1^m,\dots,x_n^m)]&\eqcirc &
R^{km}_0[ M_L(x_1^m,\dots,x_n^m)]\\
&\eqcirc &M'(x_1^m,\dots,x_n^m)
\end{eqnarray*}
where
\begin{eqnarray*}
M'(x_1,\dots,x_n)&\eqcirc &
R^{k}_0[ M_L(x_1,\dots,x_n)].
\end{eqnarray*}
A substitution of $x_1^m,\dots,x_n^m$ for $x_1,\dots,x_n$
transforms a formally zero (non-zero) polynomial
into formally zero (respectively, non-zero) polynomial, thus we get
that the graph $L$ has no coloring in $k$ or fewer colors if and
only if the polynomial
\begin{eqnarray}&
R^{k}_0[ M_L(x_1,\dots,x_n)]
\label{eq7}\end{eqnarray}
is \emph{identically equal} to the zero polynomial in the  field $\mathrm{GF}(q)$.

\vspace{3mm}
3. We will now set up a relationship between the coefficients of the polynomial
\eqref{eq7} and adjacent subclasses modulo $k$. Using the distributive property,
we can write the product \eqref{eq1}  as an algebraic sum of $2^m$ monomials
where $m$ is the number of edges of the graph $L$. There is a natural one-to-one
correspondence between these monomials and the orientations of the graph $L$: an
orientation $\overrightarrow{L}=(V,\overrightarrow{U})$
corresponds to the monomial resulting from selection in the factor
$x_i-x_j$  (where $\overrightarrow{x_ix_j}\in\overrightarrow{U}^*$)
either the first or the second summand depending on whether
$\overrightarrow{x_ix_j}\in\overrightarrow{U}$ or
$\overrightarrow{x_jx_i}\in\overrightarrow{U}$, that is, the monomial
$\delta x_1^{s_1}\cdots x_n^{s_n}$ where $s_i=s_{\overrightarrow{L}}^+(x_i)$
and $\delta=1$ or $-1$ depending on the agreement or disagreement of
the orientations $\overrightarrow{L}$ and $\overrightarrow{L}^*$.

This correspondence allows us to reformulate the definitions given above in the
following new terminology: two orientations are congruent modulo $k$ if and only if
the monomials corresponding to them transform to similar monomials under maximal
reduction according to the scheme $x^k\rightarrow 1$; two orientations
agree if and only if the monomials corresponding to them are of the same sign.
Thus the set of coefficients of the polynomial \eqref{eq7} is, up to their signs,
the set of differences of the cardinalities of adjacent subclasses modulo $k$.

\vspace{3mm}
4. We now complete the proof of the proposed criterion.
Let $L=(X,U)$ be a (non-directed) graph having two adjacent subclasses modulo $k$
with different numbers of orientations. As was shown above, this is
equivalent to the statement that the polynomial \eqref{eq7} has a non-zero coefficient.
Let $p$ be a prime dividing neither this non-zero coefficient nor the number $k$.
By the Dirichlet box principle, among the $k+1$ numbers $1,p,\dots,p^k$ there are two distinct
numbers congruent modulo $k$, and, respectively, $p^t\equiv 1\pmod{k}$ for some
positive $t$ (it is well-known that for such a $t$ we can take $\phi(k)$ where
$\phi$ is Euler's totient function but we need only the mere existence of such a $t$).
Let $q=p^t$, then the polynomial \eqref{eq7} isn't formally equal to the zero polynomial
 in the  field $\mathrm{GF}(q)$ because its coefficients belong to the prime
subfield, and $p$ doesn't divide any of the coefficients. As was shown above, this
implies that graph $L$ has a vertex coloring in at most $k$ colors.

Now let $L=(X,U)$ be a (non-directed) graph having a vertex coloring in at most $k$ colors,
let $l$ be a natural number different from $1$ and co-prime with $k$. Let $p$
be a prime factor of $l$. Let us select $q$ in the same manner as was done in the
proof of sufficiency. Then the polynomial \eqref{eq7} has a coefficient different from
zero in the field  $\mathrm{GF}(q)$; this coefficient is not a multiple of $p$ and
hence not a multiple of $l$. The equivalence class modulo $k$ corresponding to this
coefficient splits into two adjacent subclasses having cardinalities different
modulo $l$.  The necessity is proved.

\begin{center}
References
\end{center}

1. Vitaver L.M. Finding minimal vertex coloring of a graph with Boolean powers of
the incidence matrix (in Russian). \emph{Dokl. AN SSSR}, 1962, 147:4, pp.\,758-759;
\url{http://www.ams.org/mathscinet-getitem?mr=145509}.

2. Minty G.I. A theorem on $n$-coloring  the points of a linear graph.
\emph{Amer. Math. Monthly},  1962, 69:7, pp.\,623--624.

3. Zykov A.A. \emph{Theory of finite graphs} I. Novosibirsk, Nauka Publishing House,
1969.

\

\hfill \parbox{60mm}{Received by the editorial board\\
on June 18, 1974}

\end{document}